\documentclass[12pt]{article}  

\usepackage{color}
\usepackage{abstract}
\usepackage{amsmath}
\usepackage{amsfonts}
\usepackage{bbm}
\usepackage{amsthm}
\usepackage{comment}
\usepackage[nottoc]{tocbibind}
\usepackage{changepage}

\usepackage{natbib}
\usepackage[colorlinks=true,citecolor=blue,runcolor=blue,linkcolor=blue, pdfborder={0 0 0}]{hyperref}

\usepackage[hang,flushmargin]{footmisc}

\title{Exact testing with random permutations}
\author{Jesse Hemerik and Jelle Goeman\\
	Leiden University Medical Center, The Netherlands
	}

\theoremstyle{plain}
\newtheorem{theorem}{Theorem}

\newtheorem{proposition}{Proposition}

\theoremstyle{definition}
\newtheorem{definition}{Definition}

\newtheorem{condition}{Condition}

\theoremstyle{remark}

\renewcommand\footnotemark{}

\newcommand{\hp}{H_{\text{p}}}

\begin{document}
\maketitle




\begin{abstract}
\noindent When permutation methods are used in practice, often a limited number of random permutations are used to decrease the computational burden. However, most theoretical literature assumes that the whole permutation group is used, and methods based on random permutations tend to be seen as approximate. There exists a very limited amount of literature on exact testing with random permutations and only recently a thorough proof of exactness was given. In this paper we provide an alternative proof, viewing the test as a ``conditional Monte Carlo test'' as it has been called in the literature. We also provide  extensions of the result. Importantly, our results can be used to prove properties of various  multiple testing procedures based on random permutations. 

\let\thefootnote\relax\footnotetext{Address: Department of Medical Statistics and Bioinformatics, Leiden University Medical Center, Postzone S5-P, Postbus 9600, 2300RC Leiden, The Netherlands.\\ E-mail: j.b.a.hemerik@lumc.nl}

\end{abstract}



\section{Introduction}
\label{intro}

Permutation tests are nonparametric tests that are used in particular when the null hypothesis implies distributional invariance under certain transformations \citep{fisher1936coefficient, lehmann2005testing, ernst2004permutation}. 
Apart from permutations, other groups of transformations can be used, such as rotations \citep{langsrud2005rotation}. 

When the set of transformations used is not a group, a permutation test can be very conservative or anti-conservative. The first author who explicitly assumed a group structure is \cite{hoeffding1952large}. The role of the group structure has recently been emphasized \citep{southworth2009properties, goeman2010sequential}. \citet{southworth2009properties} note that in particular the set of `balanced permutations' cannot be used, since it is not a group.

Often it is computationally infeasible to use the whole group of permutations, due to its large cardinality.
In that case random permutations are used, as was first proposed by \cite{dwass1957}. Often a permutation \emph{p}-value based on random permutations is simply seen as an estimate of the permutation \emph{p}-value.

It is known that naively using random permutations instead of all possible permutations can lead to extreme anti-conservativeness \citep{phipson2010permutation}, especially when combined with multiple testing procedures.
Therefore sometimes the identity permutation, which corresponds to the original observation, is included with the random permutations \citep{ge2003resampling, lehmann2005testing}.  \cite{lehmann2005testing} (p.636) state that when the identity is added, the estimated \emph{p}-value is stochastically larger than the uniform distribution on $[0,1]$ under the null. 
\cite{phipson2010permutation} note that adding the identity can make the permutation test exact, i.e. of level $\alpha$ exactly. 
They do not mention the role of the underlying group structure.
Instead they view the permutation test as a Monte Carlo test, which is known to be exact in some situations if the original observation is added.

Referring to Monte Carlo is not sufficient, because despite being related, a Monte Carlo test is very different from a permutation test. 
Monte Carlo samples are draws from the null distribution. 
In the permutation context, the random permutations of the data are instead drawn from a conditional null distribution, i.e. the permutation distribution.
Hence the proof by \cite{phipson2010permutation} is incomplete and it remained unclear what assumptions (e.g. a group structure) are essential for the validity of random permutation tests. For example, it is unclear from Phipson and Smyth that random sampling from balanced permutations would lead to invalid tests.

In \citet{hemerik2017false} a test is given based on random transformations.
In the present paper we extend this work, investigating fundamental properties of random permutation tests. Our main focus is on the level of tests.
Other properties, e.g. consistency, do not generally hold but can be established for more specific scenarios \citep{lehmann2005testing, pesarin2015some, pesarin2013weak} by using results presented here. Our results are general and can be used to prove properties of various multiple testing methods based on random permutations, such as \citet{westfall1993resampling}, \citet{tusher2001significance},  \citet{meinshausen2005lower} and  \citet{meinshausen2006false}.
In the literature there are two approaches to proving permutation tests with fixed permutations: a conditioning-based approach \citep{pesarin2015some} and a more direct approach  \citep{hoeffding1952large, lehmann2005testing}. We will give proofs with both approaches.

The structure of the paper is as follows. In section \ref{sec1} we review known results on the level of a permutation test based on a fixed group of transformations.  The concepts and definitions from section \ref{sec1} are used throughout the paper. Testing with random permutations is covered in Section \ref{sec2}. In section \ref{secmc} permutation tests are contrasted with Monte Carlo tests. 
Estimation of \emph{p}-values is discussed in section \ref{secest}.
Exact tests and \emph{p}-values based on random transformation are given in section \ref{secrandom} and \ref{secpv2}. In section \ref{secap} some additional applications of these results are mentioned.

\section{Fixed transformations} \label{sec1}

Here we discuss tests that use the full group of transformations.

\subsection{Basic permutation test} \label{secbasic}

Let $X$ be data taking values in a sample space $\mathcal{X}$. 
Let $G$ be a finite set of transformations $g:\mathcal{X}\rightarrow\mathcal{X}$, such that $G$ is a group with respect to the operation of composition of transformations.
This means that $G$ satisfies the following three properties: $G$ contains an identity element (the map $x\mapsto x$); every element of $G$ has an inverse in $G$; for all $a_1,a_2\in G$, $a_1 \circ a_2 \in G$. 
This assumption of a group structure for $G$ is fundamental throughout the paper, since it ensures that $Gg=G$ for all $g\in G$, i.e. that the set $G$ is permutation-invariant.

Considering a general group of transformations rather than only permutations is useful, since in many practical situations the group  consists of e.g. rotations \citep{langsrud2005rotation, solari2014rotation} or maps that multiply part of the data by $-1$ (\cite{pesarin2010permutation}, pp.\ 54 and 168). 
We write $g(X)$ as $gX$. Consider any  test statistic $T:\mathcal{X}\rightarrow \mathbb{R}$. Throughout this paper we are concerned with testing the following null hypothesis of permutation-invariance.

\begin{definition}
Let $\hp$ be any null hypothesis which implies that the joint distribution of the test statistics $T(gX)$, $g\in G$, is invariant under all transformations in $G$ of $X$. That is, writing $G=\{a_1,...,a_{\#G}\}$, under $\hp$
\begin{equation} \label{pinv}
\big(T(a_1X),...,T(a_{\#G}X)\big)\,{\buildrel d \over =}\, \big(T(a_1gX),...,T(a_{\#G}gX)\big)
\end{equation}
for all $g\in G$.
\end{definition}
Note that \eqref{pinv} holds in particular when for all $g\in G$ $$X\,{\buildrel d \over =}\,gX.$$ Composite null hypotheses are usually not of the form $\hp$, but for specific scenarios, properties of tests of such hypotheses can be established using results in this paper.

The most basic permutation test rejects $\hp$ when $T(X)>T^{(k)}(X)$, where 
$$ T^{(1)}(X)\leq...\leq T^{(\#G)}(X)$$ are the sorted test statistics $T(gX)$, $g\in G$, and $k=\lceil (1-\alpha)\#G\rceil$ with $\alpha\in[0,1)$. As is known and stated in the following theorem, this test has level at most $\alpha$. 

\begin{theorem} \label{basic}
Under $\hp$, $\mathbb{P}\big\{T(X)>T^{(k)}(X)\big\}\leq \alpha.$
\end{theorem}

We now give two proofs: a conditioning-based approach and an approach without conditioning. Both approaches are more or less known.  The conditioning-based proof is similar to that in \cite{pesarin2015some} but the setting is more general. For each $x\in \mathcal{X}$, define $O_x$ to be the orbit of $x$, which is the set $\{gx:g\in G\}\subseteq \mathcal{X}.$

\begin{proof}
Let $A=\{x\in \mathcal{X}: T(x)>T^{(k)}(x)\}$ be the set of elements of the sample space that lead to rejection.
Suppose $\hp$ holds. By the group structure, $Gg=G$ for all $g\in G$. Consequently, $T^{(k)}(gX)=T^{(k)}(X)$ for all $g\in G$.
Thus 
$\#\{g\in G: gX\in A\}=$ $$ \#\{g\in G: T(gX)>T^{(k)}(gX)\}= \#\{g\in G: T(gX)>T^{(k)}(X)\}\leq\alpha\#G .$$

Endow the space of orbits with the $\sigma$-algebra that it  inherits from the $\sigma$-algebra on  $\mathcal{X}$. 
Analogously to the proof of theorem 15.2.2 in \citet{lehmann2005testing}, we obtain
$$\mathbb{P}(X\in A  \mid  O_X)= \frac{1}{\#G}\#\{ g\in G: gX\in A \}. $$
By the argument above, this is bounded by $\alpha$.
Hence $$\mathbb{P}(X\in A)=E\big\{\mathbb{P}(X\in A \mid O_X)\big\}\leq\alpha$$ 
as was to be shown. 
\end{proof}

We now state a different proof without conditioning. A similar proof can be found in \cite{hoeffding1952large} and \cite{lehmann2005testing} (p. 634).
\begin{proof}
By the group structure, $Gg=G$ for all $g\in G$. Hence $T^{(k)}(gX)=T^{(k)}(X)$ for all $g\in G$.
Let $h$ have the uniform distribution on $G$. Then under $\hp$, the rejection probability is
$$
\begin{aligned}
&\mathbb{P}\big\{T(X)>T^{(k)}(X)\big\}=\\
&\mathbb{P}\big\{T(hX)>T^{(k)}(hX)\big\}=\\
&\mathbb{P}\big\{T(hX)>T^{(k)}(X)\big\}.
\end{aligned}
$$
The first equality follows from the null hypothesis and the second equality holds since $T^{(k)}(X)=T^{(k)}(hX)$. 
Since $h$ is uniform on $G$, the above probability equals
$$\mathbb{E} \Big [(\#G)^{-1} \cdot \#\big \{g\in G: T(gX)> T^{(k)}(X) \big \} \Big ] \leq \alpha,$$ 
as was to be shown.

\end{proof}

The test of theorem \ref{basic} is not always exact.
When the data are discrete then the basic permutation test is often slightly conservative, due to a non-zero probability of tied values in $X$. 
Under the following condition, which is often satisfied for continuous data, but usually not for discrete data, the test is exact for certain values of $\alpha$. 

\begin{condition} \label{cond1}
There is a partition $\{G_1,...,G_m\}$ of $G$ with $id\in G_1$ and $\#G_1=...=\#G_m$, such that under $\hp$ with probability $1$ for all $g,g'\in G$,  $T(gX)=T(g'X)$ if and only if $g$ and $g'$ are in the same set $G_i$.
\end{condition}

\begin{proposition} \label{basicexact} 
Under condition \ref{cond1}, the test of theorem \ref{basic} is exact  if and only if  $\alpha \in \{0,1/m,...,(m-1)/m\}.$ 
\end{proposition}

The proof of this result is analogous to the proof of theorem \ref{basic}. 
As an example where condition \ref{cond1} holds, consider a randomized trial where $X\in \mathbb{R}^{2n}$ and the test statistic is
\begin{equation}
T(X)=\sum_{i=1}^{n} X_i - \sum_{i=n+1}^{2n} X_i,   \label{tst}
\end{equation} 
 where $X_1,...,X_{n}$ are cases and $X_{n+1},...,X_{2n}$ are controls and all $X_i$ are independent and identically distributed under the null. Let $$m=\binom{2n}{n}.$$ If the observations are continuous then the set of $\alpha$ for which the test is exact is  $\big\{ 0, 1/m,..., (m-1)/m \big\},$
reflecting the fact that there are $m$ equivalence classes of size $n!n!$ of permutations that always give the same test statistic.

The test of theorem \ref{basic} is often conservative when the data are discrete, since then condition \ref{cond1} is usually not satisfied.  
Moreover, in many cases, the value 0.05 is not in the set mentioned in proposition \ref{basicexact} and hence the permutation test for $\alpha=0.05$ is conservative, even if condition \ref{cond1} is satisfied.
The test can be adapted to be exact by randomizing it, i.e. by rejecting  $\hp$ with a suitable probability $a$ in the boundary case that $T(X)=T^{(k)}$ \citep{hoeffding1952large}. Here
\begin{equation} \label{eqa}
a=a(X)= \frac{\alpha \#G-M^+(X) }{M^0(X)},
\end{equation}
where $$M^+(X):=\#\{g\in G: T(gX)>T^{(k)}(X)\},$$
$$M^0(X):=\#\{g\in G: T(gX)=T^{(k)}(X)\}.$$

This adaptation has the advantage that it is always exact.  Even if condition \ref{cond1} is satisfied, the adaptation can be useful to guarantee that the level of the test is exactly the nominal level $\alpha$.  On the other hand, this test is less reproducible than the test of theorem \ref{basic}, since its outcome may depend on a random decision. Which test is to be preferred, would depend on the context.

When the set $G$ is not a group, the test can be highly anti-conservative or conservative. For example, the set of \emph{balanced permutations} is a subset of the set of all permutations, but is not a subgroup. These permutations have been used in various papers since they can have an intuitive appeal. They are discussed in \cite{southworth2009properties}, who warn against their use  since they lead to anti-conservative tests. 
The fact that permutations have been used incorrectly illustrates that more emphasis should be put on the assumption of a group structure.

Intuitively, the reason why a group structure is needed for theorem \ref{basic} is the following. Suppose for simplicity that $H_{\text{p}}$ implies that $X \,{\buildrel d \over =}\, gX$ for all $g\in G$. The permutation test works since under $H_{\text{p}}$, for every permutation $g\in G$ the probability  $\mathbb{P}\{T(gX)>T^{(k)}(X)\}$ is the same. 
The reason is that under $H_{\text{p}}$, for every $g\in G$, the joint distribution of $gX$ and the set $GX$, i.e. of $(gX,GX)$, is the same.
Indeed, since $G=Gg$ (group structure), the set $GX$ is  a function of $gX$, namely $GX=f(gX)$, with $f$ given by $f(x)=Gx$.
Thus, for $g,g'\in G$, $(gX,GX)=(gX,f(gX)) \,{\buildrel d \over =}\, (g'X,f(g'X))= (g'X,GX)$.
 When $G$ is not a group,   the joint distribution of $gX$ and the set $GX$ is not generally independent of $g$.

\subsection{Permutation \emph{p}-values} \label{secpvfixed}

Permutation \emph{p}-values are \emph{p}-value based on permutations of the data.
Here we will discuss permutation \emph{p}-values based on the full permutation group. \emph{p}-values based on random permutations are considered in  section  \ref{secpv2}. 

It is essential to note that there is often no unique null distribution of $T(X)$, since $H_p$ often does  not  specify a unique null distribution of the data. Correspondingly,  $T^{(k)}(X)$ should not be seen as the  $(1-\alpha)$-quantile of \emph{the} null distribution.

When  a test statistic $t$ is a function (which is not random) of the data and has a unique distribution under a hypothesis $H$, then a \emph{p}-value in the strict sense,
$\mathbb{P}_H(t\geq t_{obs})$, is defined where $t_{obs}$ is the observed value of $t$.
Since under $H_p$ $T(X)$ does not always have a unique null distribution,  often there exists no \emph{p}-value in the strict sense based on this test statistic. 
However, under condition \ref{cond1} the statistic 
$$D=\#\big \{g\in G: T(gX)\geq T(X)\big \}$$
does have a unique null distribution. Thus a \emph{p}-value in the strict sense based on $-D$ is then defined. 
Denoting by $d$ the observed value of $D$, we have under $\hp$
$$\mathbb{P} (-D\geq -d)=\mathbb{P}(D\leq d)=\mathbb{P}\big\{T(X)>T^{(\#G-d)}(X)\big\}=\frac{d}{\#G}.$$
This is indeed what is usually considered to be the permutation \emph{p}-value.
This equality holds under condition \ref{cond1}.
In other cases, such as when the observations are discrete, the null hypothesis often does not specify a unique null distribution of $D$. Thus there is not always a \emph{p}-value in the strict sense based on $D$.

When $H_p$ does not specify a unique null distribution of any sensible test statistic, as a resolution a `worst-case' \emph{p}-value could be defined. However sometimes better solutions are possible, e.g. the randomized \emph{p}-value $p'$ in section  \ref{secpv2}. 
In general, a \emph{p}-value in the weak sense can be considered, i.e. any random variable $p$ satisfying $\mathbb{P}(p\leq c)\leq c$ for all $c\in [0,1]$ for every distribution under the null hypothesis. For $\hp$, $D/\#G$ is always a \emph{p}-value in the weak sense.

\section{Random transformations} \label{sec2}

In section  \ref{sec2} we extend the results of the previous section to tests based on random transformations. 
Since permutation testing with random permutations is often confused with Monte Carlo testing, in Section \ref{secmc} the differences between the two are made explicit. Since random permutations are often used for estimation (rather than exact computation) of \emph{p}-values, estimation of permutation \emph{p}-values is discussed in section \ref{secest}.
Exact tests and \emph{p}-values are given in sections \ref{secrandom} and  \ref{secpv2} respectively. These two sections contain most of the novel results of the paper.

\subsection{Comparison of Monte Carlo and permutation tests} \label{secmc}

In a basic Monte Carlo experiment, the null hypothesis $H_0$ is that $X$  follows a specific distribution. A Monte Carlo test is used when there is no analytical expression for the $(1-\alpha)$-quantile of the null distribution of $T(X)$, such that the observed value of $T(X)$ cannot simply be compared to this quantile. To test $H_0$, independent realizations $X_2,...,X_w$ are drawn from the null distribution of $X$. Assume that $T(X),T(X_2),...,T(X_w)$ are continuous. 
Writing $X_1=X$, let 
$$B'=\#\{1\leq j \leq w: T(X_j)\geq T(X)\}$$ 
and let $b'$ denote its observed value.
It is easily verified that under $H_0$, $B'$ has the uniform distribution on $\{1,...,w\}$.

The Monte Carlo test rejects $H_0$ when  $T(X)>T^{(k')}$, where $k'=\lceil (1-\alpha) w\rceil$ and $T^{(1)}\leq ... \leq T^{(w)}$ are the sorted test statistics $T(X_j)$, $1\leq j \leq w$.  
Note that $T^{(k')}$ is not the exact $(1-\alpha)$-quantile of the null distribution of $T(X)$, but nevertheless the test is exact. 
The reason is that the null distribution of $B'$ is known. The test rejects $H_0$ if and only if $-B'$ exceeds the $(1-\alpha)$-quantile of its null distribution.  Equivalently, it rejects when the Monte Carlo \emph{p}-value $$\mathbb{P}_{H_0}(B'\leq b')=b'/w,$$ where $b'$ is the observed value of $B'$, is at most $\alpha$. 

The validity of a random permutation test is not as obvious. Let $g_2,...,g_w$ be random permutations from $G$. (There are various ways of sampling them, which we discuss later.) One permutation, $g_1$, is fixed to be $id\in G$, reflecting the original observation. Then, similarly to a Monte Carlo test, the permutation test rejects $\hp$ if and only if $T(X)>T^{(k')}(X)$, 
where now
$T^{(1)}\leq ... \leq T^{(w)}$ are the sorted test statistics $T(g_j X)$, $1\leq j \leq w$.

Note that contrary to the Monte Carlo sample $X_1,...,X_w$, the permutations $g_1X,...,g_wX$ are not independent  under the null. 
Thus the random permutation test is not analogous to the  Monte Carlo test. To prove the validity of the test based on random permutations, we must use that $g_1X,...,g_wX$ are independent and identically distributed conditionally on the orbit $O_X$. It is however not obvious what properties $G$ should have in order that $g_1X=X$ can be seen as a random draw from $O_X$ conditionally on $O_X$. It will be seen that it suffices that $G$ is a group.
In that case, the test can be said to be a `conditional Monte Carlo test'.


\subsection{Estimated \emph{p}-values} \label{secest}
In practice it is often computationally infeasible to calculate the permutation \emph{p}-value based on the whole permutation group, $D/ \#G$.
To work around this problem, there are two approaches in the literature. In both approaches, random permutations are used.
The first approach is \emph{calculating} (rather than estimating) a \emph{p}-value based on the random permutations. 
This is discussed in section  \ref{secpv2}.
The second approach is \emph{estimating} the \emph{p}-value $D/ \#G$, which we discuss now.

In practice the  \emph{p}-value $p=D/\#G$ is often estimated using random permutations.  The random permutations  are typically all taken to be uniform on $G$ and can be drawn with or without replacement. The estimate of $p$ is often taken to be 
$\hat{p}=B/w$, with $B$ as defined above.
This is an unbiased estimate of $p$, i.e. $E\hat{p}=p$, and usually $\lim_{w\rightarrow\infty}\hat{p}=p$.

A more conservative estimate $\tilde{p}=(B+1)/(w+1)$ is sometimes also used. This formula is discussed in section \ref{secpv2}.

Using the unbiased estimate  $\hat{p}=B/w$ can be very dangerous, as \citet{phipson2010permutation} thoroughly explain. The reason is that $\hat{p}$ is almost never stochastically larger than the uniform distribution on $[0,1]$ under $\hp$. This is immediately clear from the fact that $\hat{p}$ usually has a strictly positive probability of being zero. Consequently, if $\hp$ is rejected if $\hat{p}\leq c$ for some cut-off $c$, then the type-I error rate can be larger than $c$. Often this difference will be small for large $w$. However, when $c$ is itself small due to e.g. Bonferroni's multiple testing correction, then $\mathbb{P}(\hat{p}\leq c)$ can become many times larger than $c$ under $\hp$. This is because this probability does not converge to zero as $c\downarrow 0$ for fixed $w$. Thus, as 
\cite{phipson2010permutation} note, using $\hat{p}$ in combination with e.g. Bonferroni can lead to completely faulty inference.
Appreciable anti-conservativeness also occurs if very few (e.g.\ 25--100) random permutations are used  
(as in e.g. \cite{byrne2013monozygotic} and \cite{schimanski2013tracking}).

When possible, computing exact \emph{p}-values is always to be preferred over estimating \emph{p}-values. Exact \emph{p}-values based on random permutations are given in section \ref{secpv2}.

\subsection{Random permutation tests} \label{secrandom}
Here we discuss exact tests based on random transformations. Apart from theorem \ref{main} \citep{hemerik2017false}, the results in this section are novel.  

\citet{phipson2010permutation} also consider exact \emph{p}-values based on  random permutations.
The proofs in \citet{phipson2010permutation} are incomplete, since they do not show the role of the group structure of the set of all permutations.
\citet{lehmann2005testing} (p.636) remark without proof that if $G$ is a group, then under $\hp$ the \emph{p}-value $(B+1)/(w+1)$ is always stochastically larger than uniform on $[0,1]$, but they state no other properties.
In \citet{hemerik2017false}  for the first time a theoretical foundation is given for the random permutation test, using the group structure of the set $G$. 
Here this work is extended with additional results.

Theorem \ref{main} states that the permutation test with random permutations has level at most $\alpha$ if the identity map is added. This was remarked several times in the literature and proved in \citet{hemerik2017false}. We first define the vector of random transformations.

\begin{definition} \label{defrp}
Let $G'$ be the vector $(id, g_2, ..., g_w)$, where $id$ is the identity in $G $ and $g_2, ..., g_w$ are random elements from $G$. Write $g_1=id$. The transformations can be drawn either with or without replacement: the statements in this paper hold for both cases. If we draw $g_2,...g_w$ \emph{without} replacement, then we take them to be uniformly distributed on $G\setminus \{id\}$, otherwise uniform on $G$. In the former case, $w\leq \#G$.
\end{definition}

\begin{theorem} \label{main}
Let $G'$ be as in Definition \ref{defrp}.
Let $T^{(1)}(X,G')\leq...\leq T^{(w)}(X,G')$ be the ordered test statistics $T(g_j X)$, $1\leq j \leq w$.
Let $\alpha \in [0,1)$ and recall that $k'=\lceil (1-\alpha)w \rceil$.

Reject $H_p$ when $T(X,G')>T^{(k')}(X,G')$.
Then the rejection probability under $\hp$ is at most $\alpha$.
\end{theorem}

A proof of  theorem \ref{main} is in \citet{hemerik2017false} and we recall it here.

\begin{proof}
From the group structure of $G$, it follows that for all $1\leq j \leq w$, $G'g_j^{-1}$ and $G'$ have the same distribution, if we disregard the order of the elements.
Let $j$ have the uniform distribution on $\{1,...,w\}$ and write $h=g_j$.
Under $\hp$,
$$
\begin{aligned}
&\mathbb{P}\big\{ {T}(X)> T^{(k')}(X,G')\big\}=\\
&\mathbb{P}\big\{ {T}(X)> T^{(k')}(X,G'h^{-1})\big\}=\\
&\mathbb{P}\big\{ {T}(hX)> T^{(k')}(hX,G'h^{-1})\big\}.\\
\end{aligned}
$$
Since $(G'h^{-1})(hX)=G'(h^{-1}hX)$, the above equals
$$
\begin{aligned}
&\mathbb{P}\big\{ {T}(hX)> T^{(k')}(h^{-1}hX,G')\big\}=\\
&\mathbb{P}\big\{ {T}(hX)> T^{(k')}(X,G')\big\}.\\
\end{aligned} $$
Since $h=g_j$ with $j$ uniform, this equals
$$\mathbb{E} \Big [w^{-1} \#\big \{1\leq j \leq w: T^{(j)}(X,G')> T^{(k')}(X,G') \big \} \Big ] \leq \alpha,$$ 
as was to be shown.
\end{proof}

We now prove theorem \ref{main} with a conditioning-based approach, viewing the test as a ``conditional Monte Carlo" test as it has been called in the literature.

\begin{proof} We prove the result for the case of drawing with replacement. The proof for drawing without replacement is analogous.
Note that $(X,G')$ takes values in $\mathcal{X}\times\{id\}\times G^{w-1}$.
Let $A\subset \mathcal{X}\times\{id\}\times G^{w-1}$ be such that the test rejects if and only if $(X,G')\in A$.

Endow the space of orbits with the $\sigma$-algebra that it inherits from the $\sigma$-algebra on $\mathcal{X}$. 
Suppose $\hp$ holds. 
Assume that almost surely $O_X$ contains $\#G$ distinct elements. In case not, the proof is analogous.
Analogously to the proof of theorem 15.2.2 in \cite{lehmann2005testing}, we obtain
\begin{equation} \label{eq:cp}
\mathbb{P}\big \{(X,G')\in A\mid O_X\big \}= \frac{\#\big(O_X\times \{id\}\times G^{w-1}  \big)\cap A}{\# O_X\times \{id\}\times G^{w-1}}.        
\end{equation}

We now argue that this is at most $\alpha$. Fix $X$. Let $\tilde{X}$ have the uniform distribution on $O_X$. It follows from the group structure of $G$ that the entries of $G' \tilde{X}$ are just independent uniform  draws from $O_X$.  Thus from the Monte Carlo testing principle it follows that $\mathbb{P}\{(\tilde{X},G')\in A\big\}\leq \alpha.$ Since $(\tilde{X},G')$ was uniformly distributed on $O_X\times \{id\}\times G^{w-1}$, it follows that \eqref{eq:cp} is at most $\alpha$.  We conclude that $$\mathbb{P}\big \{(X,G')\in A\big \}   =   E\Big [\mathbb{P}\big \{(X,G')\in A\mid O_X\big \} \Big ]\leq \alpha,$$
as was to be shown. 
\end{proof}

Theorem \ref{main} implies that $(B+1)/(w+1)$ is always a \emph{p}-value in the weak sense if all random permutations (including $g_1$)  are uniform draws with replacement from $G$ or without replacement from $G\setminus \{g_1\}.$ Under more specific assumptions, theorem \ref{main} can be extended to certain composite null hypotheses. 
Proposition \ref{condition} states that under condition \ref{cond1} and suitable sampling, the test with random permutations is exact. 
The formula in Section \ref{secpv2} for the \emph{p}-value under sampling without replacement is equivalent to the last part of this result.

\begin{proposition} \label{condition}
Suppose condition \ref{cond1} holds. Let $h_1\in G_1 ,...,h_m\in G_m$.
Then the result of theorem \ref{main} still holds if $g_2,...,g_w$ are drawn with replacement from $\{h_1,...,h_m\}$ or without replacement from 
$\{h_2,...,h_m\}.$
Moreover, in the latter case, the test of theorem \ref{main} is exact for all $\alpha\in\{0/w,1/w, ... ,(w-1)/w\}.$
\end{proposition}

\begin{proof}
We consider the case that $g_2,...,g_w$ are drawn without replacement from $\{h_2,...,h_m\}$ and show that the test is exact for $\alpha\in \{0/w,...,(w-1)/w\}$. Write $G'=(g_1,...,g_w).$ 
Let $h$ have the uniform distribution on $\{g_1,...,g_w\}$.
For each $g\in G$ let $i(g)\in\{1,...,m\}$ be such that $g\in G_{i(g)}$.
Suppose $\hp$ holds. From the group structure of $G$ it follows that the sets $\big \{i(g_1),...,i(g_w)\big \}$ and $\big \{i(g_1h^{-1}),...,i(g_wh^{-1})\big \}$ have the same distribution. Consequently
$$
\begin{aligned}
&\mathbb{P}\big\{{T}(X)> T^{(k')}(X,G')\big\}=\\
&\mathbb{P}\big\{{T}(X)> T^{(k')}(X,G'h^{-1})\big\}.\\
\end{aligned}
$$
As in the above proof of theorem \ref{main} we find that this equals $\mathbb{P}\big\{ {T}(hX)> T^{(k')}(X,G')\big \}$. 

Since $\alpha\in \{0/w,...,(w-1)/w\}$ and $T(g_1X),...,T(g_wX)$ are distinct, it holds with probability one that
$$\#\big \{1\leq j \leq w: T(g_jX)>T^{(k')}\big \}=\alpha w.$$
Since $h$ is uniform, it follows that $\mathbb{P}\big\{ T(hX)> T^{(k')}(X,G')\big\}=\alpha.$    
\end{proof}

Using this result it can be shown that specific tests with random permutations are unbiased. 
The test of theorem \ref{main} can be slightly conservative if $\alpha$ is not chosen suitably or due to the possibility of ties. Recall that the same holds for the basic permutation test that uses all transformations in $G$. The adaptation by \citeauthor{hoeffding1952large} at \eqref{eqa} then guarantees exactness. The following is a generalization of \citeauthor{hoeffding1952large}'s result to random transformations.

\begin{proposition} \label{randomized}
Consider the setting of theorem \ref{main}. Let
\begin{equation} \label{eqa2}
a=a(X,G')= \frac{w\alpha-M^+(X,G') }{M^0(X,G')},
\end{equation}
where $$M^+(X,G'):=\#\{1\leq j \leq w: T(g_jX)>T^{(k')}(X,G')\},$$
$$M^0(X,G'):=\#\{1\leq j \leq w: T(g_jX)=T^{(k')}(X,G')\}.$$
Reject if $T(X)>T^{(k')}(X,G')$ and reject with probability $a$ if $T(X)=T^{(k')}(X,G')$. Then the rejection probability is exactly $\alpha$ under $\hp$.
\end{proposition}

\begin{proof}
Assume $\hp$ holds. Note that
$$\mathbb{P}(\text{reject})= E\big \{ \mathbbm{1}_{\{T(X)> T^{(k')}(X,G')\}}+a(X,G')\mathbbm{1}_{\{T(X)= T^{(k')}(X,G')\}}\big \}.$$
Write $M^+=M^+(X,G')$ and $M^0=M^0(X,G')$. Analogously to the first four steps of the first proof of theorem \ref{main}, it follows for $h$ as defined there that the above equals
$$
\begin{aligned}
&E\big \{\mathbbm{1}_{\{T(hX)> T^{(k')}(X,G')\}}+a(X,G')\mathbbm{1}_{\{T(hX)= T^{(k')}(X,G')\}}\big \}=\\
&E\big \{\mathbbm{1}_{\{T(hX)> T^{(k')}(X,G')\}} \big \} +E \big \{a(X,G')\mathbbm{1}_{\{T(hX)= T^{(k')}(X,G')\}}\big \}=\\
&E\big \{  M^+ w^{-1}\big\}     +  E\Big [E \Big \{\frac{w\alpha-M^+ }{M^0}\mathbbm{1}_{\{T(hX)= T^{(k')}(X,G')\}} \mid M^+,M^0\Big \}\Big ]=\\
&E\big \{  M^+ w^{-1}\big\}     +  E\Big [ \frac{w\alpha-M^+ }{M^0} E\Big \{\mathbbm{1}_{\{T(hX)= T^{(k')}(X,G')\}}\mid M^+,M^0\Big \} \Big ]=\\
&E\big \{  M^+ w^{-1}\big\}     +  E\Big [ \frac{w\alpha-M^+}{M^0} M^0 w^{-1} \Big ]=\alpha,\\
\end{aligned}
$$ 
as was to be shown. 
\end{proof}

The test of proposition \ref{randomized} entails a randomized decision: in case $T(X)=T^{(k')}$, the test randomly rejects with probability $a$. This is in itself not objectionable, since the test is randomized anyway due to the random transformations. 
Note that in the situation of proposition \ref{condition} under drawing without replacement the test is already exact, such that proposition \ref{randomized} is not needed to obtain an exact test.

In theorem \ref{main}, the requirement of using the whole group is replaced by suitable random sampling from the group. Interestingly, the following sampling scheme is also possible. Let $G^*\subseteq G$ be any finite subset of $G$, where we now allow $G$ to be an infinite group as well. Write $k^*=\lceil (1-\alpha)\#G^*\rceil$.
Let $h$ be uniformly distributed on $G^*$ and independent. Reject $\hp$ if and only if $$T(X)>T^{(k^*)}(X,G^*h^{-1}),$$
i.e. if $T(X)$ exceeds the $(1-\alpha)$-quantile of the values $T(gh^{-1})$, $g\in G^*$.  This is a randomized rejection rule, since it depends on $h$, which is randomly drawn each time the test is executed. The rejection probability is at most  $\alpha$, which follows from an argument analogous to the last five steps of the first proof of theorem \ref{main}. Note that if $G^*$ is a group itself, then $G^*h^{-1}=G^*$ and this test becomes non-random, coinciding with the basic permutation test. Thus it is a generalization thereof. 
This result allows using a permutation test when $G$ is an infinite group of transformations, from which it may not be obvious how to sample uniformly. 
One simply uses any finite subset $G^*$ of the infinite group.

\subsection{\emph{p}-values based on random transformations} \label{secpv2}

\citet{phipson2010permutation} give formulas for \emph{p}-values, when permutations are randomly drawn.  
Here we provide the required assumptions and proofs, which follow from section \ref{secrandom}. We then provide some additional results. 

Write
\begin{equation} \label{eq:defB}
{B}=\# \big\{1\leq j \leq w: T({g_j}X)\geq T(X)\big\},
\end{equation}
where ${g}_1,...,{g}_w$ are random permutations with distribution to be specified. Let $b$ be the observed value of $B$. Under condition \ref{cond1}, Phipson and Smyth's \emph{p}-values are exactly equal to $\mathbb{P}_{\hp}(-B\geq -b)$. 
Under condition \ref{cond1}, if ${g}_1,...,{g}_w$ are drawn such that they are from distinct elements $G_i$ of the partition and not from $G_1$, the \emph{p}-value $\mathbb{P}_{\hp}(-B\geq -b)$ is exactly 
$$\frac{b+1}{w+1}.$$
The validity of this formula follows from proposition \ref{condition}.
For the case that permutations are drawn with replacement, where ${g}_1,...,{g}_w$ are independent and uniform on $G$,  Phipson and Smyth also provide a formula for $\mathbb{P}_{\hp}(-B\geq -b)$, under condition \ref{cond1}. 

The formula $(B+1)/(w+1)$ simplifies to the formula $B/w$ if the identity map is added to the random permutations.
It follows that  the permutation test based on random permutations becomes exact for certain $\alpha$ if the identity is added.
Note that this only holds if condition \ref{cond1} is satisfied and all permutations are from distinct equivalence classes $G_i$.

We now state some additional results that follow from section \ref{secrandom}.
Corresponding to the randomized test of proposition \ref{randomized}, a randomized \emph{p}-value can be defined as follows.
The advantage of this \emph{p}-value is that it is always uniform on $[0,1]$ under $\hp$ without requirement of additional assumptions, and it is easy to compute.
Consider the randomized test of proposition \ref{randomized} (hence with $G'$ as in definition \ref{defrp}).
Suppose without loss of generality that when $T(X)=T^{(k')}$, the test rejects if and only if $a>u$, where $u$ is uniform on $[0,1]$ and independent. Define the randomized \emph{p}-value by
$${p}'=\frac{\#\{1\leq j \leq w: T(g_jX)> T(X)\}}{w}+u \frac{\#\{1\leq j \leq w: T(g_jX)=T(X)   \}}{w}.$$
This \emph{p}-value has the property that ${p}'\leq \alpha$ if and only if the randomized test rejects. This implies in particular that ${p}'$ is exactly uniform on $[0,1]$ under $\hp$. 
The fact that $p'$ is randomized is in itself not objectionable, since it is randomized anyway due to the random transformations.

A simple upper bound to $p'$ is
$$
\frac{\#\{1\leq j \leq w: T(g_jX)\geq T(X)\}}{w},
$$
a \emph{p}-value in the weak sense,
which translates to $(B+1)/(w+1)$ when  ${g}_1,...,{g}_w$ are for example all independent uniform draws from $G$.
It is not exactly uniform on $[0,1]$ under $\hp$. However, when $w$ is large and there are few ties among the test statistics, it tends to closely approximate ${p}'$, so that it may be used for simplicity.

\section{Applications} \label{secap}
We briefly mention some applications where our results are particularly useful. We have considered data $X$ that lie in an arbitrary space $\mathcal{X}$ and an arbitrary group of transformations $G$. 
For example, we allow $X$ to be a vector of functions, which is the type of data investigated by functional data analysis (FDA) \citep{cuevas2014partial, goia2016introduction}. 
\citet{cox2008pointwise} consider permutation testing with such functional data. To formulate an exact random permutation test in such a setting, the present paper is useful. 

In \citet{hemerik2017false}, properties are proven of the popular method SAM \citep[``Significance Analysis of Microarrays",][]{tusher2001significance}. This is a permutation-based multiple testing method which provides an estimate of the false discovery proportion, the fraction of false positives among the rejected hypotheses. Using  theorem \ref{main},   \citet{hemerik2017false} showed for the first time how a confidence interval can be constructed around this estimate.

In a basic permutation test, the observed statistic $T(X)$ is compared to $T^{(k)}\in \mathbb{R}$, a quantile of the permutation distribution.
The permutation-based multiple testing method by \citet{meinshausen2006false}, which provides simultaneous confidence bounds for the false discovery proportion, also constructs a quantile based on the permutation distribution. There, however, $l\in \mathbb{N}$ hypotheses and hence $l$ statistics $T_1(X),....,T_l(X)$, are considered. (They consider \emph{p}-values as test statistics.) Correspondingly, the quantile which \citeauthor{meinshausen2006false} constructs is $l$-dimensional. It turns out that the crucial step of the proof (the second last line of the proof, p. 231) relies on the principle behind the basic permutation test. The present article can be used to make this method exact. (For example, in \citet{meinshausen2006false}, $id$ should be added to the random permutations.)

In \citet{goeman2011multiple}, it is suggested to combine the method by \citet{meinshausen2006false} with closed testing, which leads to a very computationally intensive method. Hence preferably only a limited number of permutations (e.g. 100) would be used. The present paper allows using such a limited number of transformations, while still obtaining an exact method.

\section*{Discussion}

This paper proves properties of tests with random permutations in a very general setting. Properties such as unbiasedness of tests of composite null hypotheses and consistency do not hold in general but may be proved for more specific scenarios. 
For fixed permutations, there are many results regarding such properties \citep{hoeffding1952large, lehmann2005testing, pesarin2010permutation, pesarin2013weak} which may be extended to random permutations.

Aside from the permutation test, there are many multiple testing methods which employ permutations, some of which are mentioned in section \ref{secap}. Another example is \citet{westfall1993resampling}. These methods are precisely based on the principle behind the permutation test. This paper can provide better insight into these procedures, when random permutations are used.

\section*{Acknowledgments}
We thank Aldo Solari and Vincent van der Noort for their valuable suggestions.

\bibliographystyle{myplainnat}
\bibliography{referenties}

\begin{thebibliography}{25}
\providecommand{\natexlab}[1]{#1}
\providecommand{\url}[1]{\texttt{#1}}
\expandafter\ifx\csname urlstyle\endcsname\relax
  \providecommand{\doi}[1]{doi: #1}\else
  \providecommand{\doi}{doi: \begingroup \urlstyle{rm}\Url}\fi

\bibitem[Byrne et~al.(2013)Byrne, Carrillo-Roa, Henders, Bowdler, McRae, Heath,
  Martin, Montgomery, Krause, and Wray]{byrne2013monozygotic}
Byrne, E., Carrillo-Roa, T., Henders, A., Bowdler, L., McRae, A., Heath, A.,
  Martin, N., Montgomery, G., Krause, L., and Wray, N.
\newblock Monozygotic twins affected with major depressive disorder have
  greater variance in methylation than their unaffected co-twin.
\newblock \emph{Translational psychiatry}, 3\penalty0 (6):\penalty0 e269, 2013.

\bibitem[Cox and Lee(2008)]{cox2008pointwise}
Cox, D.~D. and Lee, J.~S.
\newblock Pointwise testing with functional data using the westfall--young
  randomization method.
\newblock \emph{Biometrika}, 95\penalty0 (3):\penalty0 621--634, 2008.

\bibitem[Cuevas(2014)]{cuevas2014partial}
Cuevas, A.
\newblock A partial overview of the theory of statistics with functional data.
\newblock \emph{Journal of Statistical Planning and Inference}, 147:\penalty0
  1--23, 2014.

\bibitem[Dwass(1957)]{dwass1957}
Dwass, M.
\newblock Modified randomization tests for nonparametric hypotheses.
\newblock \emph{The Annals of Mathematical Statistics}, 28:\penalty0 181--187,
  1957.

\bibitem[Ernst et~al.(2004)]{ernst2004permutation}
Ernst, M.~D. et~al.
\newblock Permutation methods: a basis for exact inference.
\newblock \emph{Statistical Science}, 19\penalty0 (4):\penalty0 676--685, 2004.

\bibitem[Fisher(1936)]{fisher1936coefficient}
Fisher, R.~A.
\newblock ``the coefficient of racial likeness" and the future of craniometry.
\newblock \emph{Journal of the Anthropological Institute of Great Britain and
  Ireland}, 66:\penalty0 57--63, 1936.

\bibitem[Ge et~al.(2003)Ge, Dudoit, and Speed]{ge2003resampling}
Ge, Y., Dudoit, S., and Speed, T.~P.
\newblock Resampling-based multiple testing for microarray data analysis.
\newblock \emph{Test}, 12\penalty0 (1):\penalty0 1--77, 2003.

\bibitem[Goeman and Solari(2010)]{goeman2010sequential}
Goeman, J.~J. and Solari, A.
\newblock The sequential rejection principle of familywise error control.
\newblock \emph{The Annals of Statistics}, 38:\penalty0 3782--3810, 2010.

\bibitem[Goeman and Solari(2011)]{goeman2011multiple}
Goeman, J.~J. and Solari, A.
\newblock Multiple testing for exploratory research.
\newblock \emph{Statistical Science}, 26\penalty0 (4):\penalty0 584--597, 2011.

\bibitem[Goia and Vieu(2016)]{goia2016introduction}
Goia, A. and Vieu, P.
\newblock An introduction to recent advances in high/infinite dimensional
  statistics, 2016.

\bibitem[Hemerik and Goeman(2017)]{hemerik2017false}
Hemerik, J. and Goeman, J.~J.
\newblock False discovery proportion estimation by permutations: confidence for
  significance analysis of microarrays.
\newblock \emph{Journal of the Royal Statistical Society: Series B (Statistical
  Methodology)}, 2017.

\bibitem[Hoeffding(1952)]{hoeffding1952large}
Hoeffding, W.
\newblock The large-sample power of tests based on permutations of
  observations.
\newblock \emph{The Annals of Mathematical Statistics}, 23:\penalty0 169--192,
  1952.

\bibitem[Langsrud(2005)]{langsrud2005rotation}
Langsrud, {\O}.
\newblock Rotation tests.
\newblock \emph{Statistics and computing}, 15\penalty0 (1):\penalty0 53--60,
  2005.

\bibitem[Lehmann and Romano(2005)]{lehmann2005testing}
Lehmann, E.~L. and Romano, J.~P.
\newblock \emph{Testing statistical hypotheses}.
\newblock Springer Science \& Business Media, 2005.

\bibitem[Meinshausen(2006)]{meinshausen2006false}
Meinshausen, N.
\newblock False discovery control for multiple tests of association under
  general dependence.
\newblock \emph{Scandinavian Journal of Statistics}, 33\penalty0 (2):\penalty0
  227--237, 2006.

\bibitem[Meinshausen and B{\"u}hlmann(2005)]{meinshausen2005lower}
Meinshausen, N. and B{\"u}hlmann, P.
\newblock Lower bounds for the number of false null hypotheses for multiple
  testing of associations under general dependence structures.
\newblock \emph{Biometrika}, 92\penalty0 (4):\penalty0 893--907, 2005.

\bibitem[Pesarin(2015)]{pesarin2015some}
Pesarin, F.
\newblock Some elementary theory of permutation tests.
\newblock \emph{Communications in Statistics-Theory and Methods}, 44\penalty0
  (22):\penalty0 4880--4892, 2015.

\bibitem[Pesarin and Salmaso(2013)]{pesarin2013weak}
Pesarin, F. and Salmaso, L.
\newblock On the weak consistency of permutation tests.
\newblock \emph{Communications in Statistics-Simulation and Computation},
  42\penalty0 (6):\penalty0 1368--1379, 2013.

\bibitem[Pesarin and Salmaso(2010)]{pesarin2010permutation}
Pesarin, F. and Salmaso, L.
\newblock \emph{Permutation tests for complex data: theory, applications and
  software}.
\newblock John Wiley \& Sons, 2010.

\bibitem[Phipson and Smyth(2010)]{phipson2010permutation}
Phipson, B. and Smyth, G.~K.
\newblock Permutation p-values should never be zero: calculating exact p-values
  when permutations are randomly drawn.
\newblock \emph{Statistical applications in genetics and molecular biology},
  9\penalty0 (1):\penalty0 39, 2010.

\bibitem[Schimanski et~al.(2013)Schimanski, Lipa, and
  Barnes]{schimanski2013tracking}
Schimanski, L.~A., Lipa, P., and Barnes, C.~A.
\newblock Tracking the course of hippocampal representations during learning:
  when is the map required?
\newblock \emph{The Journal of Neuroscience}, 33\penalty0 (7):\penalty0
  3094--3106, 2013.

\bibitem[Solari et~al.(2014)Solari, Finos, and Goeman]{solari2014rotation}
Solari, A., Finos, L., and Goeman, J.~J.
\newblock Rotation-based multiple testing in the multivariate linear model.
\newblock \emph{Biometrics}, 70\penalty0 (4):\penalty0 954--961, 2014.

\bibitem[Southworth et~al.(2009)Southworth, Kim, and
  Owen]{southworth2009properties}
Southworth, L.~K., Kim, S.~K., and Owen, A.~B.
\newblock Properties of balanced permutations.
\newblock \emph{Journal of Computational Biology}, 16\penalty0 (4):\penalty0
  625--638, 2009.

\bibitem[Tusher et~al.(2001)Tusher, Tibshirani, and
  Chu]{tusher2001significance}
Tusher, V.~G., Tibshirani, R., and Chu, G.
\newblock Significance analysis of microarrays applied to the ionizing
  radiation response.
\newblock \emph{Proceedings of the National Academy of Sciences}, 98\penalty0
  (9):\penalty0 5116--5121, 2001.

\bibitem[Westfall and Young(1993)]{westfall1993resampling}
Westfall, P.~H. and Young, S.~S.
\newblock \emph{Resampling-based multiple testing: Examples and methods for
  p-value adjustment}.
\newblock John Wiley \& Sons, 1993.

\end{thebibliography}

\end{document}